\documentclass[10pt]{amsart}
\usepackage{amsmath,amscd,amsthm,amssymb,enumerate, url, eucal,     rotating}
\usepackage{fullpage}

\newcommand{\negpar}[1][-1em]{%
  \ifvmode\else\par\fi
  {\parindent=#1\leavevmode}\ignorespaces
}

\swapnumbers
\theoremstyle{plain}
\newtheorem{theorem}{Theorem}[section]

\newtheorem{corollary}[theorem]{Corollary}
\newtheorem{proposition}[theorem]{Proposition}
\newtheorem{scl}[theorem]{The general cutvertex lemma}
\newtheorem{WsCorollary}[theorem]{Whitehead's cutvertex lemma}

\theoremstyle{definition}

\newtheorem{notes on whiteheads article}[theorem]{Notes on Whitehead's article}
\newtheorem{more notes on whiteheads article}[theorem]{More notes on Whitehead's article}

\newtheorem{subroutine}[theorem]{Subroutine}

\newtheorem{algorithmcv}[theorem]{Whitehead's cutvertex algorithm}

\newtheorem{definitions}[theorem]{Definitions}

\newtheorem{examples}[theorem]{Examples}

\newtheorem{notation}[theorem]{Notation}

\newtheorem{remarks}[theorem]{Remarks}

\newcommand{\EE}{\operatorname{E}}
\newcommand{\VV}{\operatorname{V}\mkern-2mu}

\DeclareMathOperator{\support}{\hskip-4pt-support}
\DeclareMathOperator{\length}{\hskip-4pt-length}
\DeclareMathOperator{\turns}{\hskip-4pt-turns}
\DeclareMathOperator{\White}{\mathbb{W}}

\def \bigast {\operatornamewithlimits{\text{\LARGE $\ast$}}}

\def\d1{\discretionary{-}{}{-}}
\def\e1{\discretionary{--}{}{--}}

\def\and1{\hskip2pt{\scriptsize\&}\hskip1pt}

\renewcommand{\phi}{\varphi}
\renewcommand{\le}{\leqslant}
\renewcommand{\ge}{\geqslant}

\begin{document}
\title{On Whitehead's first free-group algorithm,\\
cutvertices, and free-product fac\-tor\-iza\-tions}

\author{Warren Dicks\hskip1pt*}
\address{Departament de Matem{\`a}tiques,  \newline  Universitat Aut{\`o}noma de Barcelona, \newline
 E-08193 Bellaterra (Barcelona),  SPAIN \newline \vspace{-5pt}\null}
\email{dicks@mat.uab.cat \newline\normalfont {\emph {Home page}}:\,\,\,\,\url{http://mat.uab.cat/~dicks/}}
\thanks{*Research supported by MINECO (Spain) through projects numbered MTM2014-53644-P and MTM2017-83487-P}
\date{}


\begin{abstract}      Let $F$ be any  finite-rank free group, and   $R$ be
 any finite subset of \mbox{$\bigl\{g, [g]: g \in F{-}\{1\}\bigr\}$}, where
\mbox{$[g]:= \{fgf^{-1}:f\in F\}$}.
By an $R$-allocating $F$-factorization  we mean a set  \mbox{$\mathcal{H}$}
 of nontrivial  subgroups   of~$F$ such that
\mbox{$\operatornamewithlimits{\bigast}\limits_{H \in \mathcal{H}} H   =F$} \vspace{-.3mm}and
\mbox{$R \subseteq \bigl\{h, [h] : h \in H, H\in \mathcal{H}\bigr\}$}.
We show that  Whitehead's (fast) cutvertex algorithm  inputs the pair \mbox{$(F,R)$} and  outputs
a maximum-size $R$-allocating $F$-factorization.     Richard Stong showed this in  the
case where \mbox{$R \subseteq F$} or \mbox{$R\subseteq \{[g] : g \in F\}$}, thereby
unifying and   generalizing a collection of results obtained by   Berge,  Bestvina,
  Lyon,  Shenitzer,   Stallings,  Starr, and   White\-head.
   Our   proof is based on the interaction between  two normal forms for the  elements~of~$F$, rather than
the algebraic topology of handlebodies,  trees, or graph folding.
\end{abstract}

\maketitle

   \section{Outline}\label{sec:outline}

Throughout this article,   let $F$ be any finite-rank free group,
$X$ be
any $F\mkern-1mu$-basis, and
$R$~be any finite set that consists of nontrivial $F\mkern-1mu$-elements and nontrivial $F\mkern-1mu$-classes.
(An \textit{$F\mkern-1mu$-basis} is a free-gener\-ating set for~$F$,
 an \mbox{\textit{$F\mkern-3mu$-element}} is an element of $F$, and a (nontrivial)
 \mbox{\textit{$F\mkern-2mu$-class}}
is the conjugacy class of a (nontrivial) \mbox{$F\mkern-1mu$-element}.)

For \mbox{$f,g \in F$} and \mbox{$B\subseteq F$}, we write \mbox{$\null^f\mkern-2mug:= fgf^{-1}$},
\mbox{$\null^F\mkern-3mug:= \{\null^h\mkern-2mug: h \in F\}$}, and
\mbox{$\null^{\{F\}}\mkern-3muB:= \{\null^F\mkern-2mub: b \in B\}$}.

Suppose that $\mathcal{H}$ is a  multiset of subgroups of $F$ such that the  induced map
 \mbox{$\operatornamewithlimits{\bigast}\limits_{H \in \mathcal{H}} H  \to F$} is an isomorphism.
Here, we say that $\mathcal{H}$ is an  \textit{$F\mkern-1mu$-fac\-tor\-iza\-tion}, and
call the elements of  $\mathcal{H}$  \textit{the factors}.  Sometimes
we also say that the expression  \mbox{$\operatornamewithlimits{\bigast}\limits_{H \in \mathcal{H}} H$}
is an   $F\mkern-1mu$-fac\-tor\-iza\-tion. We say that the   $F\mkern-1mu$-fac\-tor\-iza\-tion \mbox{$\mathcal{H}$}
is   \textit{$R$-allo\-cating} if  each element of  \mbox{$\mathcal{H}$}  is nontrivial and
\mbox{$R \subseteq \{h, \null^F\mkern-2muh: h \in H, H\in \mathcal{H}\}$}.
   If, moreover,
 no proper free\d1product   refinement of \mbox{$\operatornamewithlimits{\bigast}\limits_{H \in \mathcal{H}} H$}  is
 $R$-allo\-cating, we say that the  $R$-allo\-cating $F\mkern-1mu$-fac\-tor\-iza\-tion \mbox{$\mathcal{H}$}~is
 \textit{atomic}.
Since $F$~has finite rank,  atomic $R$-allocating $F\mkern-1mu$-fac\-tor\-iza\-tions exist,
and we want to be able to find one as quickly as possible.

In $\S$\ref{sec:back}, we review the earlier results on this topic, starting with Whitehead's cutvertex lemma.

In $\S$\ref{sec:factorizations}, we show that
 atomic $R$-allo\-cating $F\mkern-1mu$-fac\-tor\-iza\-tions are all as similar to each other as may  reasonably be expected.
When \mbox{$R\subseteq F$},  each
atomic $R$-allo\-cating $F\mkern-1mu$-fac\-tor\-iza\-tion  gives
the  unique  inclusion-smallest free-product factor of $F$ which  includes~$R$;
this is the only information it gives when \mbox{$\vert R \vert = 1$}, where
 we may as well assume that \mbox{$R\subseteq F$}.

 In $\S$\ref{sec:RcupX}, we give a  careful treatment of  concepts
introduced by White\-head\mbox{(1936-01,\,$\S$2)}.
We denote by  \mbox{$\mathcal{P}(X;R)$} the  finest partition of~$X$ that respects the $X\mkern-3mu$-support
 of each element   of~$R$, which means that the $F\mkern-1mu$-factor\-iza\-tion
 \mbox{$\bigast\limits_{Y \mkern-3mu\in \mathcal{P}(\mkern-3muX;R)} \langle \mkern3mu Y  \mkern1mu \rangle$} \vspace{-.5mm}
is  $R$-allo\-cating; this appeared in work of  Hoare\and1Karrass\and1Solitar(1971,\,$\S$2) with different notation and
terminology.
Using \mbox{$\mathcal{P}(X;R)$} and Whitehead's graphs, we  define
\textit{$R$-cut\-ver\-tex\d1free} $F\mkern-1mu$-bases.
We then present \textit{Whitehead's cut\-vertex algorithm},~$\S$\ref{alg:ca} below,
which inputs the pair~\mbox{$(X,R)$} and outputs
an $R$-cutvertex\d1free $F\mkern-1mu$-basis.  (By an \textit{algorithm}  we
 mean  a procedure with choices whose possible outputs  have some specified property.)

 In $\S$\ref{sec:main}, we  present \textit{the   general cutvertex\,\,lemma},~$\S$\ref{lem:scl} below,   which says
\begin{equation}\label{eq:GCL}
\text{if $X$ is $R$-cut\-vertex\d1free, then the  $R$-allo\-cating $F\mkern-1mu$-factor\-iza\-tion
 \mbox{$\bigast\limits_{Y \mkern-3mu\in \mathcal{P}(\mkern-3muX;R)} \langle \mkern3mu Y  \mkern1mu \rangle$} is
atomic.}\phantom{-----?--}
\end{equation}

 \vspace{-2mm}

In summary, we show that Whitehead's  (fast) cutvertex algorithm inputs the pair~\mbox{$(X,R)$} and
outputs an atomic $R$-al\-locating $F\mkern-1mu$-fac\-tor\-iza\-tion.

\section{Chronology of proofs of cases  of the   general cutvertex lemma}\label{sec:back}

$\bullet$ Whitehead(1936-01) proved the case of~\eqref{eq:GCL} where $R$ is a subset of
\mbox{$B$} or  \mbox{$\null^{\{\mkern-2muF\mkern1mu\}}\mkern-2muB$}
for some $F$-basis~$B$;
in detail, he used  the algebraic topology  of a certain\vspace{.6mm} three\d1manifold
to prove that if  $X$~is $R$-cutvertex\d1free here, then  $R$ is a subset of
\mbox{$X\, {\cup} \, X^{-1}$} or
 \mbox{$\null^{\{\mkern-2muF\mkern1mu\}}\mkern-2mu(X\, {\cup} \, X^{-1})$} respectively.
This is called \textit{Whitehead's cutvertex lemma},  \mbox{$\S$\ref{lem:wcl} below}.
Gersten(1984,\,Example) announced  that   graph\d1theoretic machinery he had developed  could be used to prove
 the \mbox{$R \subseteq X\, {\cup} \, X^{-1}$} result,
and Hoare(1988,\,Theorem~3)  provided such a proof.

Put together,  Whitehead's cutvertex algorithm and cutvertex lemma  constituted
the first-ever \textit{sub-basis algorithm}, by which we mean
  an algorithm which extends
a given finite subset of $F$ to an $F\mkern-1mu$-basis or determines that that is not possible,
and analogously for a given finite set of $F\mkern-1mu$-classes.

$\bullet$ Shenitzer(1955,\,\,Corollary) used another result of  White\-head(1936-10,\,\,Theorem~3)
to prove    the  case of~\eqref{eq:GCL}  where \mbox{$\vert R \vert = 1$} and
$X$~is  $R$-min\-i\-miz\-ing; the latter concept is defined in $\S$\ref{mnowa} below.

$\bullet$ Lyon(1980,\,\,Theorem 1)  developed  Shenitzer's method to prove    the   case of~\eqref{eq:GCL}
where $R$ is a subset of \mbox{$F$}~or~\mbox{$\null^{\{\mkern-2muF\mkern1mu\}}\mkern-2muF$}
and $X$~is  $R$-min\-i\-miz\-ing.

$\bullet$   Starr(1992) gave cutvertex arguments which, in the form distilled  by
 Wu(1996,\,\,\hskip-.7pt$\S$1),  and in the light of a result of Lyon(1980,\,\,Theorem\,2),
prove what is the case of~\eqref{eq:GCL} where $R$~is  the  set of $F\mkern-1mu$-classes  determined
by  a finite set  of disjoint simple closed curves
on the boundary of a handlebody which has~$F$ as fundamental group. It was Stallings(1999,\,$\S$3)
who realized that   Whitehead's cutvertex algorithm was being given
a completely new application here.

$\bullet$   John Berge, in a 1993 preprint,   proved
the case of~\eqref{eq:GCL}  where  \mbox{$\vert R \vert = 1$},   by using the algebraic topology of
Whitehead's three\d1mani\-fold; see
Stallings(1999, Corollary\,2.5).  Nata\hskip.3pt{\v s}\hskip.1pta Macura kindly informed me that
 Mladen Best\-vina independently proved the same case, by analyzing infinite paths in a Cayley tree; see
 Martin(1995,\,\,Theorem\,49).

$\bullet$ Stong(1997,\,\,Theorem\,3) proved the case of~\eqref{eq:GCL}   where
\mbox{$R \subseteq \null^{\{\mkern-2muF\mkern1mu\}}\mkern-2muF$}, by analyzing
bi-infinite paths in a Cayley tree.
Independently, Stallings(1999,\,\,Theorem\,2.4) proved the  same case,   by using  the algebraic topology of
Whitehead's three\d1mani\-fold.  An elegant graph-theoretic folding  proof was given by Wilton(2018, Lemma~2.10)
and, independently, by Heusener\and1Weidmann(2019, $\S$3).

$\bullet$ Stong(1997,\,\,Theorem\,10) proved the case of~\eqref{eq:GCL}   where \mbox{$R \subseteq F$},
by using the algebraic topology of a handle\-body. A Bass-Serre-theoretic two-tree proof was given by Dicks(2014, $\S$2).

\section{Atomic $R$-allo\-cating $F\mkern-1mu$-fac\-tor\-iza\-tions}\label{sec:factorizations}

 Recall that $F$ is a finite-rank free group  and $R$ is a  finite subset of \mbox{$\bigl\{g, \null^F\mkern-3mug : g \in F{-}\{1\}\bigr\}$}.

\begin{definitions} Artin(1926) gave the normal form for an element  of a free product  of groups,
as presented by Serre(1977, I.1.2.1).  This may be used to prove
that if $\mathcal{H}$  is any $F\mkern-1mu$-fac\-tor\-iza\-tion, then, for any \mbox{$H_1, H_2 \in \mathcal{H}$} and \mbox{$g_1, g_2 \in F$}, if
\mbox{$\null^{g_1}\mkern-3muH_{1} \,{\cap}\, \null^{g_2}\mkern-3muH_{2}  \ne \{1\}$} in $F$, then
\mbox{$g_1H_{1} = g_2H_{2}$} in $F$,  and   \mbox{$H_1 = H_2$} in \mbox{$\mathcal{H}$} and in~$F$.
 This implication may also be viewed as a consequence of the result of Serre(1977, I.5.3.12) that if $\mathcal{H}\ne \emptyset$,
then the disjoint union of the left $F\mkern-1mu$-sets \mbox{$F/H$},  \mbox{$H\in \mathcal{H}$},
  is the vertex-set of a left $F\mkern-1mu$-tree with   trivial edge stabilizers.

Recall that an  $F\mkern-1mu$-fac\-tor\-iza\-tion \mbox{$\mathcal{H}$}
is said to be $R$-allo\-cating  if  each element of  \mbox{$\mathcal{H}$} is nontrivial,
each   $F\mkern-2mu$-element  \mbox{$r \in R$} is an element of
 some~\mbox{$H \in \mathcal{H}$},
and   each  $F\mkern-2mu$-class \mbox{$r \in R$}  contains an element of some~\mbox{$H \in \mathcal{H}$},
which implies that \mbox{$r\,{\cap}\,H$} is an $H$-class;
 in each case, we now see that the element $H$ of \mbox{$\mathcal{H}$} is uniquely determined by~$r$,
and we shall say that \textit{$r$ is allo\-cated to
the factor $H$}.
For each \mbox{$H \in \mathcal{H}$},
we write \mbox{$R_{\vert H,\mathcal{H}}$}
to denote the set  of  elements of~\mbox{$R$} which are allo\-ca\-ted to~$H$,
 sometimes viewed  as a set of   $H$-elements and $H$-classes.
Notice that \mbox{$\{R_{\vert H,\mathcal{H}} : H \in \mathcal{H}\}-\{\emptyset\}$} is a partition of~$R$.

If \mbox{$F = \{1\}$}, then \mbox{$R = \emptyset$}, and
$\emptyset$ is the unique $R$-allo\-cating $F\mkern-1mu$-fac\-tor\-iza\-tion.
 If \mbox{$F \ne \{1\}$},  then   \mbox{$\{F\}$}
is an $R$-allo\-cating $F\mkern-1mu$-fac\-tor\-iza\-tion; if it is the only one,
then we say that $F$ is an \textit{$R$-atom}.
Thus, an  $R$-allo\-cating $F\mkern-1mu$-fac\-tor\-iza\-tion
\mbox{$ \mathcal{H}$} is  atomic
if and only  each \mbox{$H \in \mathcal{H}$} is an \mbox{$R_{\vert H,\mathcal{H}}$}-atom. \qed
\end{definitions}

 \begin{proposition}\label{prop:unique} All the atomic $R$-allo\-cating $F\mkern-1mu$-fac\-tor\-iza\-tions
have the same number of factors.  They all
induce the same partition  of~$R$.  For each \mbox{$r \in R$},
they all have the same $F\mkern-2mu$-conjugacy orbit of the
factor to which \mbox{$r$} is allo\-cated, and if $r$ is an $F\mkern-3mu$-element, they
all have the same factor to which  $r$    is allo\-cated.  All the  factors  with
no elements of $R$   allo\-cated to them are free subgroups of rank one.
\end{proposition}

 \begin{proof}  Consider any two atomic  $R$-allo\-cating $F\mkern-1mu$-fac\-tor\-iza\-tions
 \mbox{$ \mathcal{H}$} and   \mbox{$  \mathcal{K}$}, and  any \mbox{$H \in \mathcal{H}$}.
The  subgroup   theorem of  Kurosch(1934) gives an $H$-factorization
\mbox{$H_0 \ast
\bigast\limits_{K \in \mathcal{K}}\,\bigast\limits_{a \in A_{H,K}} (H  \cap \null^a\mkern-3muK)$}
where $H_0$ is a  free  group and, for each   \mbox{$K \in \mathcal{K}$},
 $A_{H,K}$ is a certain subset of $F$ such that \mbox{$1 \in A_{H,K}$} and the map
\mbox{$A_{H,K} \to H  \backslash F /K $}, \mbox{$a \mapsto H {\cdot}a{\cdot}K $}, is
bijective; see, for example,  Serre(1977, I.5.5.14).  Kurosch's theorem
may be used to prove the result of  Nielsen(1921) that $H$ is a free group.

Consider any \mbox{$r \in R_{\vert H,\mathcal{H}}$}.  There exists a (necessarily unique) \mbox{$K  \in \mathcal{K}$}
such that  if $r$ is an $F\mkern-1mu$-element, then
\mbox{$r \in H$}  and~\mbox{$r \in K$}, while if $r$ is an $F\mkern-1mu$-class,
then $r$ contains an $H$-element~$h$   and a $K$-element $k$.
In the former event,  \mbox{$r \in H   \cap \null^1\mkern-1muK $}.
In the latter event, since $r$ is an $F\mkern-1mu$-class, \mbox{$h =  \null^g\mkern-2muk$} for some $g \in F$,
and  then \mbox{$g = h'{\cdot}a{\cdot}k\mkern1mu'$} for some
 \mbox{$(h',a,k\mkern1mu') \in H  \times A_{H,K} \times K$}, and then  $r$ contains
 \mbox{$ \null^{h'^{-1}}\mkern-3muh  = \null^{a{\cdot}k'}\mkern-3mu k   \in H  \cap \null^a\mkern-3muK $}.
Thus, we obtain an $R_{\vert H,\mathcal{H}}$-allo\-cating $H$-fac\-tor\-iza\-tion
from \mbox{$H_0 \ast
\bigast\limits_{K \in \mathcal{K}}\,\bigast\limits_{a \in A_{H,K}} (H  \cap \null^a\mkern-3muK)$}
 by omitting  all the trivial  \vspace{-.5mm} factors.

Since \mbox{$  \mathcal{H} $} is an atomic $R$-allo\-cating $F\mkern-1mu$-fac\-tor\-iza\-tion,
$H$ is an \mbox{$R_{\vert H,\mathcal{H}}$}-atom, and there are three possibilities.
  If $R_{\vert H,\mathcal{H}} = \emptyset$, then  the rank of $H$ is $1$.
If $R_{\vert H,\mathcal{H}}$ contains an $H$-element, then \mbox{$H = H \cap K$}
and \mbox{$R_{\vert H,\mathcal{H}} \subseteq R_{\vert K,\mathcal{K}}$};
 by symmetry, \mbox{$K =  K \cap H  = H$} and $R_{\vert H,\mathcal{H}} = R_{\vert K,\mathcal{K}}$.
If \mbox{$R_{\vert H,\mathcal{H}}$} is  nonempty and consists of $H$-classes,
then \mbox{$H = H \cap \null^a \mkern-2mu K$}  and \mbox{$R_{\vert H,\mathcal{H}} \subseteq R_{\vert K,\mathcal{K}}$};
 by symmetry, $K$ is also included in an $F\mkern-1mu$-conjugate of $H$,
and we  see that \mbox{$H = \null^a \mkern-2mu K$} and $R_{\vert H,\mathcal{H}} = R_{\vert K,\mathcal{K}}$.
The result now follows.
\end{proof}

\begin{corollary}[Stong]\label{cor:stong} If $R$ is a subset of $B$ or \mbox{$\null^{\{\mkern-2muF\mkern1mu\}}\mkern-2muB$}
for some $F$-basis $B$, then each atomic $R$-allocating $F$-factorization equals \vspace{-1.5 mm}
\mbox{$\bigast\limits_{x  \,\in X} \langle x \rangle$} for some $F$-basis $X$ such that
$R$ is a subset of \mbox{$X$} or
\mbox{$\null^{\{\mkern-2muF\mkern1mu\}}\mkern-2muX$} respectively. \qed
\end{corollary}

\section{Whitehead's cutvertex algorithm}\label{sec:cuts} \label{sec:graphs}\label{sec:RcupX}

 Recall that $F$ is a finite-rank free group, $X$ is an $F\mkern-1mu$-basis,
 and $R$ is a finite   subset of \mbox{$\bigl\{g, \null^F\mkern-3mug : g \in F{-}\{1\}\bigr\}$}.

\begin{notation}  We write \mbox{$X^{\pm 1}:=X\, {\cup} \, X^{-1}$} and
\mbox{$X^{0, \pm 1}:=  \{1\}  \, {\cup} \,X^{\pm 1}$}.

If $r$ is any   $F\mkern-1mu$-element (resp.\,nontrivial $F\mkern-1mu$-class), then by an
\textit{\mbox{$X^{\pm 1}$}-word for $r$} we mean any  finite sequence \mbox{$(x_1,x_2,\ldots,x_n)$}
of elements of~\mbox{$X^{\pm 1}$} such that  $r$  equals (resp.\,contains) the product \mbox{$x_1x_2\cdots x_n$};
 here, we  set \mbox{$x_0:=1$} (resp.\,\mbox{$x_0:=x_n$}) and
 \mbox{$x_{n+1}:=1$} (resp.\,\mbox{$x_{n+1}:=x_1$}).  There exists some \mbox{$X^{\pm 1}$}-word
\mbox{$(x_1,x_2,\ldots,x_n)$} for~$r$ such that  \mbox{$x_i^{-1} \ne x_{i+1}$} for each
\mbox{$i \in \{0,1,2,\ldots,n\}$}, and
 such a word is  unique  (resp.\,unique up to cyclic permutation), which allows us to  define
$$
X\length(r) {:=} n,  \,\,\,\,
X\turns(r) {: =}  \bigl\{(x_0\hskip-5pt^{-1}\hskip-3pt, x_1), (x_1\hskip-5pt^{-1}\hskip-3pt, x_2), \ldots,
(x_n\hskip-5pt^{-1}\hskip-3pt, x_{n+1}) \bigr\}, \,\,\,\,
X\support(r) {:=} \{x_1,x_2,\ldots, x_n\}^{\pm 1} {\cap\,} X. \vspace{-.1 mm}
$$
We call each element of \mbox{$X\turns(r)$} an \textit{\mbox{$X\mkern-3mu$-turn}    of~$r$}.
If $r$ is the trivial $F\mkern-1mu$-class,  we define
\mbox{$X\length(r) {:=}0$}, \mbox{$X\turns(r):= \emptyset$}, and
\mbox{$X\support(r) {:=}\emptyset$}.
 It is not difficult to see  that if $r$ is any $F\mkern-1mu$-class, then, for each \mbox{$g \in r$},
 \mbox{$X\turns(r) \subseteq  X\turns(g^2)$}.

Let $R'$ be any subset of \mbox{$F \cup \null^{\{\mkern-2muF\mkern1mu\}}\mkern-2muF$}.  We write\vspace{-.7mm}
$$\textstyle\text{\mbox{$X\length(R') := \sum\limits_{r \in R'} X\length(r) \in  \{\infty, 0, 1, 2, \ldots\}$} \,\,\,and\,\,\,
\mbox{$X\turns(R') := \bigcup\limits_{r \in R'} X\turns(r) \subseteq X^{0, \pm 1}{\times}\,X^{0, \pm 1}$}}.\vspace{-1.3mm}$$
We call each element  of  \mbox{$X\turns(R')$}  an  \textit{\mbox{$X\mkern-3mu$-turn}    of~$R'$}.
We denote  by \mbox{$\mathcal{P}(X;R')$}  the finest partition of~$X$
such that, for each \mbox{$r \in R'$}, some element of the partition
includes \mbox{$X\support(r)$}.
We define an operation which takes  a  set of subsets of $X$ with  two overlapping elements and  replaces
those two subsets with their union, thereby reducing the number of subsets; if we   start  with the  (finite) set
\mbox{$\{X\support(r) : r \in R'\}\,\, \cup \,\, \bigl\{ \{x\}  : x \in X  \bigr\}$}
and apply this operation as often as possible,  then we obtain~\mbox{$\mathcal{P}(X;R')$}.
\qed
\end{notation}

\begin{definitions}
For any set $V$  and any element  \mbox{$(v,w)$} of \mbox{$V {\times\,} V$}, we say that
 \mbox{$(v,w)$} \textit{meets}  each subset of $V$ which
contains $v$ or $w$;   we say also that \mbox{$(v,w)$} \textit{meets}  \mbox{$v$} and \mbox{$w$}.

By a \textit{graph} $\Gamma$,
we mean a set $V$  together with a  subset $E$ of \mbox{$V {\times\,} V$}.
Then  $V$ is called the \textit{vertex-set} of~$\Gamma$,  denoted  \mbox{$\VV \Gamma$},
and its elements are called \textit{$\Gamma\mkern-2mu$-vertices},
while $E$~is called the \textit{edge-set} of $\Gamma$,  denoted \mbox{$\EE \Gamma$},
and its elements are called \textit{$\Gamma\mkern-2mu$-edges}.
A $\Gamma$-vertex~$v_\dag$ is said to be a  \textit{cutvertex of $\Gamma$} (\textit{in the sense \vspace{-.9mm}
of Whitehead}) if   \mbox{$\VV \Gamma{-}\{v_\dag\}$} equals the  union  of two  disjoint   nonempty
subsets $V_1$ and $V_2$ such that \mbox{$V_1 \overset{\Gamma}{\leftrightarrow}\mkern-16mu\not\mkern16muV_2$},
where this symbol means  ``\,no $\Gamma\mkern-2mu$-edge  meets both $V_1$ and $V_2$".  We shall use
 corresponding depictions of related phrases.

We let  \mbox{$\White(X,R)$} denote the   graph whose vertices are those elements of \mbox{$X^{0, \pm 1}$}
which are met by some  $X\mkern-3mu$-turn of~$R$ and whose edges are the $X\mkern-3mu$-turns of~$R$.
A cutvertex of $\White(X,R)$ is called an \textit{$X\mkern-3mu$-cutvertex of $\White(X,R)$} if it lies in \mbox{$X^{\pm 1}$},
that is, it does not equal $1$.

For each \mbox{$Y\mkern-5mu\in\mkern-3mu \mathcal{P}(X;R)$},
we set \mbox{$R_{\vert Y,X}:=   \{r \in R : X\support(r) \subseteq Y\}$}; then
\mbox{$\{R_{\vert Y,X}:Y \mkern-5mu\in\mkern-3mu \mathcal{P}(X;R)\}{-}\{\emptyset\}$}
is a partition of $R$.  We sometimes view the elements of \mbox{$R_{\vert Y,X}$}  as
\mbox{$\langle\mkern3muY \rangle$}-elements and \mbox{$\langle\mkern3muY \rangle$}-classes;
here,  \mbox{$\mathcal{P}(Y;R_{\vert Y,X}) = \{Y\}$}.
We say   $X$ is   \textit{$R$-cutvertex-free}  if, for
 each \mbox{$Y\mkern-5mu\in \mkern-3mu  \mathcal{P}(X;R)$},   \mbox{$\White(Y ,R_{\vert   Y,X  })$}
has no  $Y\mkern-3mu$-cut\-ver\-tices.

 An  $F$-basis   $X'$ is  a   \textit{Whitehead neighbour  of~$X$} if
\mbox{$X' \subseteq \{1,y^{-1}\}{\cdot}X{\cdot}\{1,y\}$}
for some \mbox{$y \in (X \cap X')^{\pm 1}$}.
\qed
\end{definitions}

\begin{examples}  If \mbox{$X=\{x,y\}$} and \mbox{$R=\{x, \null^F\mkern-2muy\}$}, then
\mbox{$\EE\mkern-3mu\White(X,R) =  \bigl\{(x^{-1}, 1), (1,x), (y^{-1},y)\bigr\}$} and
\mbox{$\White(X,R)$}  has four $X\mkern-3mu$-cutvertices.  Here,
\mbox{$\mathcal{P}(X;R) = \bigl\{\{x\},\{y\}\bigr\}$}, and \mbox{$X$} is
\mbox{$R$}-cutvertex-free.

If \mbox{$X\hskip-1.7pt=\hskip-1.7pt\{x,y\}$} and \mbox{$R\hskip-1.7pt=\hskip-1.7pt\{x^2y x ^{-1}y^{-1}\}$}, then
\mbox{$\EE\mkern-2mu\White(X,R)\hskip-1.7pt=\hskip-1.7pt\bigl\{(1,x), (x,y^{-1}), (y^{-1}, x^{-1}), (x^{-1},y), (y,1), (x^{-1},x)\bigr\}$} and \mbox{$\White(X,R)$}
  has no cutvertices.  Here,
\mbox{$\mathcal{P}(X;R) =  \{X\}$}, and \mbox{$X$} is
\mbox{$R$}-cutvertex-free.  \qed
\end{examples}

\begin{remarks}
In  \mbox{$\White(X,R)$}, each  vertex is met by some  edge; each edge meets two vertices since~\mbox{$1 \not\in R$};
 \mbox{$\VV\White(X,R) = \emptyset$} if and only if  \mbox{$R = \emptyset$}; and
  \mbox{$1 \in \VV\White(X,R)$}  if and only if  $R$~contains some $F\mkern-1mu$-element.

By partitioning~$X$ manually, we   reduce   our study  to the case where
\mbox{$\mathcal{P}(X;R) = \{X\}$}.
With the exception of the famous final phrase in the statement of Whitehead's cutvertex lemma,~$\S$\ref{lem:wcl} below,
we shall work with  \mbox{$\White(X,R)$} only in the case where
\mbox{$\mathcal{P}(X;R) = \{X\}$}.  One consequence is that  connectivity is not mentioned in our
arguments.

Suppose  that  \mbox{$\mathcal{P}(X;R)= \{X\}$}. \vspace{.2mm} Here,  \mbox{$X \ne \emptyset$}.
 If    \mbox{$R = \emptyset$},
 then \mbox{$\vert X \vert =1$} and    \mbox{$ \VV\White(X,R) = \emptyset$}.
If \mbox{$R \ne \emptyset$},  then
\mbox{$X^{\pm 1} \subseteq \VV\White(X,R)\subseteq X^{0,\pm1}$}.
  If \mbox{$\White(X,R)$} \textit{has} an $X\mkern-3mu$-cut\-ver\-tex, then
 Subroutine~\ref{sub:w} below   constructs   a Whitehead neighbour $X'$
of~$X$ such that    \mbox{$X'\,\length(R)   < X\length(R)$}.  If
 \mbox{$\White(X,R)$} has \textit{no} $X\mkern-3mu$-cut\-ver\-tices,
then Theorem~\ref{thm:main} below  says   that   $F$ is an
 $R$-atom \mbox{(Gr.\,$\overset{,}{\alpha}\mkern-4mu\acute{}\mkern5mu\tau o \mu\text{-}o \varsigma$\,\,`no cut').}
\qed
\end{remarks}

White\-head,  Stong,   Stallings, and others  gave cases of the following,
using similar ideas.

\begin{subroutine}\label{sub:w}    \textit{When \mbox{$\mathcal{P}(X;R)  = \{X\}$}
and   \mbox{$\White(X,R)$} has an $X\mkern-3mu$-cut\-ver\-tex,
 the following three-step procedure outputs a Whitehead neighbour $X'$
of~$X$ such that \mbox{$X'\mkern2mu\length(R) < X\length(R)$}}.

\medskip

 \textbf{Step 1.}   \textit{We set \mbox{$\Gamma :=  \White(X,R)$}.
  We shall see that  \mbox{$X^{\pm 1} \subseteq \VV\Gamma  \subseteq X^{0,\pm1}$}.  We
 shall\vspace{-.6mm}  find a   \mbox{$y_\dag \in  X^{\pm 1}$}  and an expression of
 \mbox{$\VV\Gamma {-}\{  y_\dag\}$}  as the union of
   two  disjoint   subsets \mbox{$Y_-$}~and~\mbox{$Y_+$} such that
 \mbox{$y_\dag^{-1} \in Y_- \overset{\Gamma}{\leftrightarrow}\mkern-16mu\not\mkern16muY_+
 \overset{\Gamma}{\leftrightarrow}y_\dag$}.}

  Here,  \mbox{$X^{\pm 1} \subseteq \VV\Gamma  \subseteq X^{0,\pm1}$}
since \mbox{$\mathcal{P}(X;R) = \{X\}$}  and
 \mbox{$\VV\Gamma  \ne \emptyset$}.
Since $\Gamma$ has an $X\mkern-3mu$-cutvertex,  we may find some \mbox{$x_\dag \in X^{\pm 1}$}
and \vspace{-.5mm} express \mbox{$\VV\Gamma {-}\{x_\dag\}$}
as the union of two  disjoint nonempty subsets \mbox{$X_-$}~and~\mbox{$X_+$}
such that~\mbox{$x_\dag^{-1}\in X_- \overset{\Gamma}{\leftrightarrow}\mkern-16mu\not\mkern16muX_+$}.
If \mbox{$X_+  \overset{\mkern5mu\Gamma}{\leftrightarrow} x_\dag$},
then \vspace{-1mm} setting \mbox{$y_\dag:=x_\dag$}, \mbox{$Y_-:= X_-$}, and \mbox{$Y_+:=X_+$} gives the desired
result; thus, we may assume that \mbox{$X_+  \overset{\Gamma}{\leftrightarrow}\mkern-16mu\not\mkern16mux_\dag$}.\vspace{1mm}
Here, \mbox{$\VV\Gamma $}  equals \vspace{-1mm} the union of  two disjoint  nonempty
subsets   \mbox{$Y_-:=  X_-  \cup \{  x_\dag\}$} and~\mbox{$X_+$} such that
\mbox{$Y_-  \overset{\Gamma}{\leftrightarrow}\mkern-16mu\not\mkern16muX_+$}.
For each \mbox{$v \in X_+$},  \mbox{$v \overset{\Gamma}{\leftrightarrow} X_+{-}\{v\}$};
in particular, \mbox{$X_+\ne\{1\}$}.
Since \mbox{$\mathcal{P}(X;R)= \{X\}$}, the set
\mbox{$\{ \langle X_+ \rangle, \langle\mkern3muY_-\rangle\}$} is not
an \mbox{$R$}-allo\-cating $F\mkern-1mu$-factor\-ization induced by a partition of~$X$,
which means that there  exists  some
  \mbox{$y_\dag \in X_+$} such that \vspace{-1.3mm} \mbox{$y_\dag^{-1} \in Y_-$}.
Now \mbox{$\VV\Gamma {-}\{  y_\dag\}$}
equals the union of   two  disjoint  subsets  \mbox{$Y_-$}
and \mbox{$Y_+:= X_+ {-}\{y_\dag\}$} such that
\mbox{$ y_\dag^{-1} \in Y_- \overset{\Gamma}{\leftrightarrow}\mkern-16mu\not\mkern16muY_+
 \overset{\Gamma}{\leftrightarrow}y_\dag$}.   Step  1 is completed.  We
have no further need of the two  original hypotheses.

\medskip

 \textbf{Step 2.} \textit{We have found \mbox{$y_\dag \in  X^{\pm 1}$},  and we shall construct an
\mbox{$\ell \in \{ -1,0\}$} and a map
 \mbox{$\chi:   X^{0, \pm 1} \to \{\ell, \ell{+}1\}$}, \mbox{$v \mapsto \chi(v)$},  such that
 the following hold: \mbox{$\chi( 1) = 0$};   \mbox{$\chi( y_\dag) = \chi(  y_\dag^{-1}) = \ell$};
at least one $X\mkern-3mu$-turn  \mbox{$(v,w)$} of~$R$ is {\normalfont$\chi$-cut} in the sense that
 \mbox{$\chi(v) \ne \chi(w)$};  and  every $\chi$-cut $X\mkern-3mu$-turn  of~$R$ meets $y_\dag$.}

We   set \mbox{$\ell:= -\bigl\vert \{1\} \cap Y_+ \bigr\vert$},
and form the map  \mbox{$\chi$}  which carries
\mbox{$\{1\}$} to \mbox{$\{0\}$}, \mbox{$Y_-   \cup \{ y_\dag\}  $} to \mbox{$\{\ell\}$}, and
  \mbox{$Y_+$} to \mbox{$\{\ell{+}1\}$};
our choice of $\ell$ ensures that $\chi$ is well-defined.
The \mbox{$\chi$}-cut \mbox{$X\mkern-3mu$}-turns of~$R$
are then the \mbox{$\Gamma $}-edges   which meet   \mbox{$y_\dag$}~and~\mbox{$Y_+$}, of which there
exists  at least one.    Step 2 is  completed.  We have no further need of Step~1.

\medskip

 \textbf{Step 3.}   \textit{We shall construct a Whitehead neighbour \mbox{$X'$} of $X$
such that \mbox{$X'\,\length(R) < X\length(R)$}.}

 We  set \mbox{$x':= y_\dag^{- \chi( x)}{\cdot} x {\cdot} y_\dag^{\mkern2mu \chi(  x^{-1})}$}  for
\vspace{-.5mm} each \mbox{$x \in   X^{\pm 1}$}, and
set  \mbox{$X':= \{x' : x \in X\}$}.
Then \mbox{$(x')^{-1} = (x^{-1})'$}, \mbox{$y_\dag' = y_\dag$}, and
 \mbox{$X'$} is a Whitehead neighbour of~$X$.  We consider an arbitrary \mbox{$r \in R$}, and
let \mbox{$(x_1, x_2, \ldots,x_n)$} be a shortest possible \mbox{$X^{\pm 1}$}-word for $r$.
Then \mbox{$n \ge 1$}.
If  $r$ is an $F\mkern-1mu$-element,  we set \mbox{$x_0 := x_{n+1} := 1$};
if $r$~is an $F\mkern-1mu$-class, we set \mbox{$x_0 := x_n$} and \mbox{$x_{n+1}:= x_1$}.
For each \mbox{$i \in \{0,1,\ldots, n\}$}, \mbox{$( x_{i}^{-1},  x_{i+1})$} is an
\mbox{$X\mkern-3mu$}-turn of~$r$,  and we let  \mbox{$k_i$}~denote
the unique element of the subset \mbox{$\chi\bigl(\{  x_{i}^{-1},   x_{i+1}\}-\{ y_\dag\}\bigr)$} of
\mbox{$\{\ell,\ell+1\}$}.
 For each \mbox{$i \in \{1,2,\ldots, n\}$}, we set \mbox{$g_i:=  y_\dag^{-k_{i-1}} x_i  y_\dag^{k_i}$}, and we then have
the following trichotomy:

\hskip 25ptPossibility 1: \mbox{$x_i =   y_\dag$} and  \mbox{$( x_{i-1}^{-1},   x_i)$}   is $\chi$-cut;
equivalently,  \mbox{$k_{i-1} \ne \chi(x_i)$}.

\hskip85ptHere \mbox{$k_i = \chi(x_i^{-1})$} since \mbox{$x_i^{-1} \ne y_\dag$},   and then
 \mbox{$g_i = y_\dag^{-k_{i-1}} x_i  y_\dag^{k_i}= y_\dag^{-(\ell +1)  + (1) + (\ell)} =   1$}.

\hskip25ptPossibility  2: \mbox{$x_i =  y_\dag^{-1}$} and  \mbox{$(x_{i}^{-1},   x_{i+1})$} is  $\chi$-cut;
equivalently,  \mbox{$k_{i} \ne \chi(x_{i}^{-1})$}.

\hskip85ptHere \mbox{$k_{i-1}   = \chi(x_{i})$} since   \mbox{$x_i \ne y_\dag$},
and then  \mbox{$g_i  =  y_\dag^{-k_{i-1}} x_i   y_\dag^{k_{i}}= y_\dag^{-(\ell) + (- 1) + (\ell+1)}  =  1$}.

\hskip25ptPossibility  3:   \mbox{$k_{i-1}  = \chi(  x_i)$} and  \mbox{$k_{i}  =   \chi(  x_i^{-1})$}.

\hskip85ptHere
  \mbox{$g_i = y_\dag^{-k_{i-1}} x_i  y_\dag^{k_i}
= y_\dag^{-\chi(    x_i)} x_i  y_\dag^{\mkern2mu\chi(  x_i^{-1})}
= x_i\mkern-1mu'$}. \vspace{1mm}

\noindent Since \mbox{$n=X\length(r)$}, we see that
\mbox{$X'\,\length(g_1g_2\cdots g_n) \le   X\length(r)$}  and that if
 some \mbox{$X\mkern-3mu$}-turn  of $r$  is $\chi$-cut,
then  \mbox{$X'\,\length(g_1g_2\cdots g_n) <  X\length(r)$}.
Now \nopagebreak \vspace{-3pt}
$$\textstyle g_1g_2\cdots g_n = (y_\dag^{-k_{0}} x_1  y_\dag^{k_1} )(y_\dag^{-k_{1}} x_2  y_\dag^{k_2} ) \cdots
 (y_\dag^{-k_{n-1}} x_n  y_\dag^{k_n} ) =
y_\dag^{-k_0} x_1x_2 \cdots x_n  y_\dag^{k_n}.\vspace{-4pt}$$
If  $r$ is an $F\mkern-1mu$-element,  then  \mbox{$k_0  = \chi(  x_0^{-1}) = \chi(1) = 0$} and
\mbox{$k_n  = \chi( x_{n +1})= \chi(1) = 0$};
here, \mbox{$g_1g_2\cdots g_n = r$}.  If  $r$~is an $F\mkern-1mu$-class, then
  \mbox{$k_0  =  k_n$};
here, \mbox{$g_1g_2\cdots g_n \in r$}. Thus,
\mbox{$X'\,\length(r) \le X'\,\length(g_1g_2\cdots g_n)$}.
Since at least one \mbox{$X\mkern-3mu$}-turn  of $R$ is $\chi$-cut, we see that
\mbox{$X'\,\length(R) < X\length(R)$}.  \qed
\end{subroutine}

\begin{algorithmcv}\label{alg:ca} Given \mbox{$(X,R)$}, we ask if there exists
some \mbox{$Y\mkern-5mu \in\mkern-3mu  \mathcal{P}(X;R)$} such that
 \mbox{$\White(Y,R_{\vert   Y,X  })$}
has a  $Y\mkern-3mu$-cut\-vertex.   If yes, then
 Subroutine~\ref{sub:w} outputs
   a   Whitehead neighbour  \mbox{$Y\mkern1mu'$}  of~\mbox{$Y$}
such that   \mbox{$Y\mkern1mu'\,\length(R_{\vert Y,X})
 <Y\,\length(R_{\vert Y,X})$}, and we  start  anew with $X$~replaced with
its Whitehead neighbour  \mbox{$X':= \bigl(X{-} Y \bigr) \cup Y\mkern1mu'$},
for which  \mbox{$X'\,\length(R) < X\length(R)$}.  If no,
we output~$X$,  and then stop.  This algorithm eventually outputs an  $R$-cutver\-tex\d1free $F$-basis, and then stops. \qed
\end{algorithmcv}

\begin{notes on whiteheads article}
In the types of graphs constructed by  Whitehead(1936-01,$\S2$),
each edge is given a multiplicity
and each copy of the edge is  divided into at least three edges by adding  new vertices.
It is important for his arguments that cutvertices are added to his   versions
of, for example,   \mbox{$\White\bigl(\{x\}, \{\null^F\mkern-3mux\}\bigr)$} and \mbox{$\White\bigl(\{x,y\}, \{x,y\}\bigr)$}.

  In the one sentence where he dealt with connected subgraphs,
Whitehead overlooked one case, and we work with \mbox{$\mathcal{P}(X;R)$} largely to handle that case.
Stong(1997) and, independently, Stallings(1999)  handled  it  by working with connected subgraphs
(where their term ``cut vertex"  conforms to standard usage).
With any of these straightforward  rectifications,  Whitehead's argument gives a  valid
  algorithm.
\qed
\end{notes on whiteheads article}

\section{The general  cutvertex lemma}\label{sec:main}

  Whitehead,  Stong,   Stallings, and others proved cases of the following,
using mainly the algebraic topology of handlebodies; we use the interaction between two
  normal forms.\vspace{-1mm}

\begin{theorem}\label{thm:main}   If \hskip2pt \mbox{$\mathcal{P}(X;R) = \{X\}$} and $F$ is  not an $R$-atom, then
  \mbox{$\White(X,R)$} has an $X\mkern-3mu$-cut\-vertex. \vspace{-1mm}
\end{theorem}

\begin{proof}[Proof  {\normalfont  (essentially following Dicks(2014, $\S$2))}.]
Set \mbox{$\Gamma :=  \White(X,R)$}.
As \mbox{$\mathcal{P}(X;R) = \{X\}$},  we have  \mbox{$F \ne \{1\}$} and,  also,
 if \mbox{$\vert X \vert \ne 1$}  then \mbox{$\VV\Gamma\supseteq  X^{\pm 1}$}.
As $F$ is not an $R$-atom and \mbox{$F \ne \{1\}$}, there exists
some $R$-allo\-cating  $F\mkern-1mu$-factoriza\-tion \mbox{$H_1{\ast} H_2$}.  Thus,
\mbox{$\VV\Gamma \supseteq X^{\pm 1}$}  since \mbox{$\vert X \vert \ne 1$}.

 In the case where \mbox{$\VV\Gamma =X^{\pm 1}$},  $R$  consists of $F\mkern-1mu$-classes, hence
\mbox{$\null^a\mkern-3muR = R$}  for each \mbox{$a \in F$},  and, by  replacing
 \mbox{$\{H_1,  H_2\}$} with a suitable \mbox{$\{\null^a\mkern-3muH_1,  \null^a\mkern-3muH_2\}$},
we may assume that the shortest \mbox{$X^{\pm 1}$}-word \mbox{$(y_1, y_2,\ldots, y_{n_0})$}   for some
element of  \mbox{$(H_1 \cup H_2){-}\{1\}$}  has  \mbox{$y_{n_0} \ne y_1^{-1}$}.
Taking \mbox{$y:= y_1$} and \mbox{$y':= y_{n_0}^{-1}$}, we have \vspace{-1mm}
\begin{equation}\label{eq:or}
 \text{if\, \mbox{$\VV\Gamma \ne X^{0,\pm 1}$}, then there exist  $X\mkern-3mu$-turns  \mbox{$(1,y)$}
and \mbox{$(y',1)$}
 of \mbox{$ H_1\cup H_2$} with \mbox{$1 \ne y  \ne  y' \ne 1$}.}\phantom{-----}\vspace{-1mm}
\end{equation}

For \mbox{$f \in F = H_1{\ast}H_2$},  there exists a
unique finite sequence  \mbox{$(h_1,h_2,\ldots,h_m)$}
of nontrivial elements of \mbox{$H_1 {\cup}\hskip2.2ptH_2$} such that \mbox{$h_1h_2{\cdots}\,h_m =  f$} and
 neither $H_1$ nor $H_2$
contains two consecutive terms of the sequence.  Here, we set \vspace{-4mm}
$$
\delta(f)  :=  \Bigl\{h_j^{-1}h_{j-1}^{-1}{\hskip-1pt\cdot\hskip-2pt\cdot\hskip-2pt\cdot\hskip1pt} h_1^{-1} :
j \in \{1,2,\ldots,m{-}1\}\Bigr\}
\qquad \text{and} \qquad   \chi(f) := \begin{cases}
1 &\text{if $f \ne 1$  and $h_1 \in H_1$},\\
0 &\text{otherwise.}\end{cases}  \vspace{-3mm}
$$
Clearly,  \vspace{-.5mm}
\begin{equation}\label{eq:clearly}
\text{$\delta(f)$ is finite and, also,   $\delta(f) = \emptyset$  if and only if \mbox{$f \in H_1 {\cup}\hskip2ptH_2$.}}
\phantom{-------------------}\vspace{-.5mm}
\end{equation}
We  have defined a  map  \mbox{$\chi: F \to \{0,1\}$}, \mbox{$g \mapsto \chi(g)$}, and
it is not difficult to use induction  on $m$ to prove that
\vspace{-.5mm}
\begin{equation}\label{eq}
\delta(f)  =\hskip3pt\bigl\{g \in F :  g \ne 1,\,\, gf \ne 1,\,\,\text{and}\,\, \chi(g) \ne \chi(gf) \, \bigr\}.
\phantom{----------------------}
\end{equation}

Since \mbox{$\mathcal{P}(X;R) = \{X\}$}, the \mbox{$R$}-allo\-cating  $F\mkern-1mu$-fac\-tor\-iza\-tion
\mbox{$H_1{\ast}H_2$}   is not induced by a partition of~$X$;
hence, there exists some  \mbox{$x \in X {-}(H_1{\cup}\,H_2)$}, and then
 \mbox{$\delta(x)\ne \emptyset$} by~\eqref{eq:clearly}.
Also, \mbox{$X^{\pm 1}$} is finite, and,   for each \mbox{$x \in X^{\pm 1}$},
 \mbox{$ \delta(x)$}~is finite by~\eqref{eq:clearly}.
Thus there exists some pair \mbox{$(\hat x, \hat g)$} such that
 \mbox{$\hat x \in X^{\pm 1}$}, \mbox{$\hat g \in \delta(\hat x)$}, and, subject to those two constraints,
 \mbox{$X\length(\hat g \hat x)$} is as large as possible.
By~\eqref{eq}, \mbox{$\hat g \in \delta(\hat x)$} means \mbox{$\hat g \ne 1$}, \mbox{$\hat g \hat x \ne 1$} and
\mbox{$\chi(\hat g) \ne \chi(\hat g \hat x)$}; by~\eqref{eq}, this means
 \mbox{$\hat g \hat x \in \delta(\hat x^{-1})$}.  Thus
 \mbox{$X\length(\hat g \hat x) >   X\length(\hat g)$} by the  maximality  of \mbox{$X\length(\hat g \hat x)$}.

Let  \mbox{$x_\dag$} denote the element of \mbox{$X^{\pm 1}$} such that
\mbox{$ X\length(\hat g x_\dag) < X\length(\hat g)$}.
It suffices to show that \mbox{$x_\dag$}~is a  cutvertex of \mbox{$\Gamma $}.
 Notice that \mbox{$x_\dag\ne \hat x$}.
Now \mbox{$X^{0, \pm 1}{-} \{x_\dag\}$}  equals the union of  two dis\-joint nonempty
subsets \mbox{$X_1$}~and~\mbox{$X_2$} with \mbox{$\chi(\hat g  X_1) = \{\chi(\hat g)\}$} and
\mbox{$\chi(\hat g X_2) = \{\chi(\hat g \hat x)\}$}; here
\mbox{$ 1 \in X_1$} and \mbox{$ \hat x  \in X_2$}.

\newpage

We  next prove that \nopagebreak \vspace{-1mm}
\begin{equation}\label{eq:crucial}
\text{no  \mbox{$X\mkern-4mu$}-turn  of \mbox{$H_1 {\cup}\, H_2$} meets both
\mbox{$X_1$}~and~\mbox{$X_2$}.}\phantom{-------------------------}\vspace{-1mm}
\end{equation}
\leftskip 29pt

  Consider an arbitrary \mbox{$h \in H_1 {\cup}\, H_2$}.
Let   \mbox{$(x_1, x_2, \ldots,x_n)$}  be the shortest \mbox{$X^{\pm1}$}-word for $h$.
Set  \mbox{$x_0 := x_{n+1} := 1$}.  Each \mbox{$X\mkern-3mu$}-turn of $h$ equals \mbox{$( x_{i}^{-1},  x_{i+1})$}
for some  \mbox{$i \in \{0, 1,\ldots,  n\}$}, and to prove~\eqref{eq:crucial} it suffices to show that
\mbox{$( x_{i}^{-1},  x_{i+1})$} does not meet both  \mbox{$X_1$} and \mbox{$ X_2$}.
We may assume that   \mbox{$x_\dag \not\in  \{x_{i}^{-1}, x_{i+1}\}$}, and it then suffices to show that
\mbox{$\chi(\hat g x_{i}^{-1}) = \chi(\hat g x_{i+1})$}.
Set \mbox{$h' := x_i^{-1} x_{i-1}^{-1} \cdots x_1^{-1}$} and \mbox{$h'' := x_{i+1}x_{i+2} \cdots x_n$}.  It suffices to
 show that  \mbox{$\textstyle\chi(\hat g  x_{i}^{-1})  =  \chi(\hat g h' ) =
\chi(\hat g h'')  = \chi(\hat g  x_{i+1}).$}

We first show that \mbox{$\hat g h'' \ne 1$} and
\mbox{$\chi(\hat g h'')
= \chi(\hat g  x_{i+1} )$}.
If    \mbox{$i = n$}, then  \mbox{$h'' =
  x_{n+1}x_{n+2} \cdots x_n = 1$} while  \mbox{$ x_{i+1}  = x_{n+1}  =  1$}.
The  case where \mbox{$i = n$} is now clear.
  If \mbox{$i \ne n$}, then
   \mbox{$x_{i+1}    \in     X^{\pm 1}{-}\{x_\dag\}$},
 and, hence,  \vspace{-.5mm}
   $$X\length(\hat g ) <
X\length(\hat g x_{i+1}) <
X\length(\hat g  x_{i+1}  x_{i+2}) < \cdots <
X\length(\hat g  x_{i+1}  x_{i+2}{\hskip-1pt\cdot\hskip-2pt\cdot\hskip-2pt\cdot\hskip1pt}x_n).$$
By the  maximality  of \mbox{$X\length(\hat g \hat x)$},  \mbox{$ \hat g x_{i+1} \not\in \delta(x_{i+2})$}.
By~\eqref{eq}, \mbox{$\chi(\hat g  x_{i+1}) = \chi(\hat g  x_{i+1}  x_{i+2})$}.
Continuing in this way, we  see that  \vspace{-1mm}
$$
\chi(\hat g  x_{i+1}) =
\chi(\hat g  x_{i+1}  x_{i+2}) =  \cdots = \chi(\hat g x_{i+1}
 x_{i+2}{\hskip-1pt\cdot\hskip-2pt\cdot\hskip-2pt\cdot\hskip1pt}x_n).\vspace{0mm}$$
Since  \mbox{$x_{i+1}x_{i+2} \cdots x_n  =  h'' $}, the  case where \mbox{$i \ne n$}  is now clear also.

By using $h'$ in place of $h''$ in the preceding argument, we see that \mbox{$\hat g  h'   \ne 1$}
and
 \mbox{$ \chi (\hat g h' )=  \chi (\hat g  x_{i}^{-1})$}.

By~\eqref{eq:clearly}, \mbox{$\delta(h) = \emptyset$}.  By~\eqref{eq}, \mbox{$\chi (\hat g  h' ) =  \chi(\hat g h'')$} since   \mbox{$ (\hat g  h' )  h
   =   \hat g h''$}.  Thus, \eqref{eq:crucial} holds.

\leftskip 0pt

  Notice that \mbox{$X\turns(R) \subseteq X\turns(H_1 {\cup}\, H_2)$},
since each \mbox{$r\in R$} equals-or-con\-tains  some    \mbox{$b  \in H_1{\cup}\,H_2$},  whence
  \mbox{$X\turns(r) \subseteq X\turns( \{b,b^2\})$}. Since \mbox{$\EE\Gamma  = X\turns(R)$},
it follows from~\eqref{eq:crucial}  that  no \mbox{$\Gamma $}-edge  meets both
\mbox{$X_1$}~and~\mbox{$X_2$}.  Since
\mbox{$X_2 \cup \{ x_\dag \} \subseteq X^{\pm 1} \subseteq \VV\Gamma  \subseteq X^{0, \pm 1}$},
it suffices to show that \mbox{$X_1 \cap \VV\Gamma  \ne  \emptyset$}.
We may then assume that \mbox{$\VV\Gamma \ne X^{0,\pm 1}$}, and
  here  \mbox{$\{y,y'\}{-} \{x_\dag  \}  \subseteq  X_1 \cap \VV\Gamma $}
by~\eqref{eq:or} and~\eqref{eq:crucial}.
This proves Theorem~\ref{thm:main}.
\end{proof}

\begin{scl}\label{lem:scl}
 If $F$ is a  finite-rank free group,  $R$    a finite subset of
 \mbox{$\bigl\{g,   \null^{F}\mkern-2mu  g  \hskip-1.5pt : \hskip-1.5pt  g \hskip-1pt \in \hskip-1pt F{-}\{1\}\bigr\}$},
\vspace{-1.5mm} and  \mbox{$X$}~an   $R$-cutvertex-free $F\mkern-1mu$-basis,
then \mbox{$\bigast\limits_{Y \mkern-3mu\in  \mathcal{P}(\mkern-3muX;R)} \langle \mkern3mu Y  \mkern1mu \rangle$}
is an  atomic $R$-allo\-cating $F\mkern-1mu$-fac\-tor\-iza\-tion.
\end{scl}

\begin{proof} For each \mbox{$Y \mkern-5mu\in \mkern-3mu   \mathcal{P}(X;R)$},  we know that
 \mbox{$\mathcal{P}(Y;R_{\vert Y, X})    = \{Y\}$}  and that
\mbox{$\White(Y,R_{\vert  Y,X  })$} has no $Y\mkern-3mu$-cut\-ver\-tices; hence
\mbox{$\langle\, Y \rangle$}~is an   \mbox{$R_{\vert   Y,X }\mkern2mu$}-atom by Theorem~\ref{thm:main}.
\end{proof}

\begin{WsCorollary}\label{lem:wcl}  Let $F$ be a  finite-rank free group, $B$ and $X$ be   $F$-bases,
 and $R$ be a subset of $B$ $($resp.\,\mbox{$\null^{\{\mkern-2muF\mkern1mu\}}\mkern-2muB$}$)$.
If \mbox{$X$} is $R$-cutvertex-free, then  $R$ is a subset of \mbox{$X^{\pm 1}$}
$($resp.\,\mbox{$\null^{\{\mkern-2muF\mkern1mu\}}\mkern-2mu(X^{\pm 1})$}$)$.
 Contrapositively, if   $R$ is not a subset of
\mbox{$X^{\pm 1}$} $($resp.\,\mbox{$\null^{\{\mkern-2muF\mkern1mu\}}\mkern-2mu(X^{\pm 1})$}$)$,
then \mbox{$X$} is not   $R$-cutvertex-free, whence   \mbox{$\White(X,R)$}~has an $X\mkern-3mu$-cutvertex. \vspace{-.5mm}
\end{WsCorollary}

\begin{proof}[Proof {\normalfont (Stong)}.] This follows from  Lemma~\ref{lem:scl} and Corollary~\ref{cor:stong}.
\end{proof}

\begin{more notes on whiteheads article}\label{mnowa}
Each $F\mkern-1mu$-basis has only finitely many Whitehead neighbours.
We say that  $X$ is \textit{$R$-minimizing} if
\mbox{$X\mkern1mu'\,\length(R) \ge  X\,\length(R)$}
 for each  Whitehead neighbour $X\mkern1mu'$  of~$X$.
\textit{Whitehead's  minimizing algorithm}  constructs,  by trial and error,  a
Whitehead-neighbour-choosing sequence   which  starts at $X$,
makes \mbox{$\operatorname{length}(R)$}   smaller with each step,
 and arrives at an $R$-minimizing
$F\mkern-1mu$-basis in at most  \mbox{$X\length(R)-\vert R \vert$}  steps.
 Whitehead's cutvertex algorithm,~$\S$\ref{alg:ca} above, shows that
 $R$-minimizing $F\mkern-2mu$-bases are  $R$-cutvertex-free.
The \textit{only} fact highlighted by Whitehead(1936-01)
was  that his    min\-imizing   algorithm  and his cutvertex lemma, when put together,
  constituted    a sub-basis algorithm which could be applied to the study of three\d1manifolds.
 His  cutvertex    algorithm produces  a  faster sub-basis algorithm that is somewhat overlooked,
perhaps because White\-head(1936-10,\,\,Theorem~3) later
 proved an extremely useful result about $R$-mini\-mizing $F\mkern-3mu$-bases.   (Although
  stated   for the case where
 \mbox{$R \subseteq F$} or \mbox{$R \subseteq  \null^{\{\mkern-2muF\mkern1mu\}}\mkern-2muF$},
the result was proved for the case where
 \mbox{$R \subseteq F \cup \null^{\{\mkern-2muF\mkern1mu\}}\mkern-2muF$}; see, for example, Dicks(2017).)
If matches had been invented   after  the cigarette lighter,
they would have been  the sensation of the twentieth century,
and if  Whitehead had written (1936-01) after  (1936-10), his cutvertex algorithm  would have
become the gold standard of sub-basis algorithms.  \nopagebreak

The marvellous \textit{Quote Investigator}
gives several published variants of the assertion about matches;   the earliest one
is attributed to Charles Norris.  \, \url{https://quoteinvestigator.com/2017/12/06/matches/}  \qed
\end{more notes on whiteheads article}

\bigskip

 \pagebreak

\centerline{\sc References}

 \bigskip

\leftskip 18pt \parindent -18pt %in foregoing:  \showthe \leftskip 0pt \parindent 12pt

E.\,Artin:\,\,Das freie Produkte von Gruppen.\,\,pp.\,361--364 in \mbox{Felix} Klein:\,\,Vorlesungen
\"uber h\"ohere Geo\-metrie.  Dritte Auflage.  Bearbeitet und herausgegeben von \mbox{W\hskip-1pt.}\,Blaschke.
Grundlehren Math.\,\,Wiss.\,\,\textbf{22}. Springer, Berlin, viii{+}405 pages (1926).

%John\,Berge:\,\,Full-rank closure.  Preprint\hskip1pt\footnote{\hskip1ptCredited by Stallings(1999, Corollary 2.5).} (1993).

% Mladen\,Bestvina:\,\,Whitehead graphs and decomposition spaces. Personal communication to Reiner
% Martin\hskip1pt\footnote{\hskip1pt Credited by Martin(1995, Theorem 49).}  (1995).

%Bruce\,Chandler and Wilhelm\,Magnus:\,\,The history of combinatorial group theory.\,\,A case study in the
%history of ideas.\,\,\,Stud.\,Hist.\,Math.\,and Phys.\,Sci.\,\textbf{9}. Springer, New York,  viii{+}234 pages (1982).

% Cooke,\,\,Alistair:\,\,Talk about America 1951--1968.  The Bodley Head, London (1968)

Warren\,Dicks:\,\,On free-group algorithms that sandwich a subgroup between free-product factors.
J.\,\,Group Theory \textbf{17}, 13--28   (2014).

Warren\,Dicks:\,\,A graph-theoretic proof for Whitehead's second free-group algorithm.\,\,14 pages\,\,(2017).\\
\url{https://arxiv.org/abs/1706.09679}

%Warren\,Dicks and M.\,J.\,Dunwoody:\,\,Groups acting on graphs.
%Cambridge Studies in Advanced Mathematics  \textbf{17}, CUP, Cambridge xvi+283 pages\,\,(1989).
%Errata at \url{http://mat.uab.cat/~dicks/DDerr.html}

%Walther\,Dyck:\,\,Gruppentheoretische Studien.\,\,\,Math.\,Ann.\,\textbf{20}, 1--44 (1882).

%Gersten,\,\,S.\,\,M.:\,\,On fixed points of certain automorphisms of free groups.
%Proc.\,\,London   Math.\,\,Soc.\,\,(3)
%\textbf{48}, 72--90  (1984a).

%Ralph\,H.\,Fox:\,\,Free differential calculus I.\,\,Derivation in the free group ring.
%Ann.\,Math.\,\textbf{57}\,(2), 547--560  (1953).

S.\,M.\,Gersten:\,\,On Whitehead's algorithm. Bull.\,\,Amer.\,\,Math.\,\,Soc.\,\,\textbf{10},  281--284    (1984).

 Michael\,Heusener\,and\,Richard\,Weidmann:\,\,A remark on Whitehead's cut-vertex lemma.
J.\,\,Group Theory \textbf{22}, 15--21   (2019).

A.\,H.\,M.\,Hoare:\,\,On automorphisms of free groups\,\,I.\,\,\,J.\,\,London Math.\,\,Soc.\,\,(2) \textbf{38},
277--285  (1988).

A.\,Howard\,M.\,Hoare,\,Abraham Karrass,\,and\,Donald\,Solitar:\,\,Subgroups of finite index of Fuchsian groups. Math.\,\,Z. \textbf{120},
289--298 (1971).

Alexander\,Kurosch:\,\,Die Untergruppen der freien Produkte von beliebigen Gruppen.\,Math.\,Ann.\,\textbf{109}, 647--660
  (1934).

Herbert\,C.\,Lyon:\,\,Incompressible surfaces in the boundary of a handlebody\,--\,an algorithm.
Canadian J.\,\,of Math.\,\,\textbf{32}, 590--595 (1980).

Reiner\,Martin:\,\,Non-uniquely ergodic foliations of thin type, measured currents and
auto\-morphisms of free groups. PhD thesis, %
%MR2693216, %
UCLA, xi{+}87 pages (1995).
\newline\url{https://search.proquest.com/docview/304185823}

%Otto\,Schreier:\,\,Die Untergruppen der freien Gruppen.
%Abh.\,Math.\,Univ.\,Hamburg\,\textbf{5}, 161--183 (1927).

J.\,Nielsen:\,\,Om Regning med ikke-kommutative Faktorer og dens Anvendelse i Gruppe\-te\-or\-ien.%
\,\,Mat.\,Tidsskrift B \textbf{107}, 77--94 (1921).

Jean-Pierre\,Serre:\,\,Arbres, amalgames, {${\rm SL}_{2}$}.\,\,Cours au
Coll\`ege de France r\'edig\'e avec la col\-laboration de\linebreak Hyman Bass.
Ast\'erisque \textbf{46}.\,\,Soc.\,\,Math.\,\,de\,\,France,  Paris, 189 pages  (1977).

 Abe\,Shenitzer:\,\,Decomposition of a group with a single defining relation into a free product.
 Proc.\,\,Amer. Math.\,\,Soc.\,\,\textbf{6}, 273--279 (1955).

 John\,R.\,Stallings:\,\,Whitehead graphs on handlebodies.
pp.\,317--330 in  Geometric group
theory \mbox{down under} (Canberra, 1996)
(eds.\,John Cossey, Charles F.\,Miller\,III, Walter D.\,Neu\-mann, and Michael Sha\-piro),
 Walter de Gruyter, Berlin, xii{+}333 pages  (1999).

Edith\,Nelson\,Starr:\,\,Curves in handlebodies. PhD thesis,  %
%MR2688492, %
UC Berkeley, iv{+}22 pages  (1992).
\newline\url{https://search.proquest.com/docview/303991902}

Richard\,Stong:\,\,Diskbusting elements of the free group. Math.\,\,Res.\,\,Lett.\,\,\textbf{4}, 201--210 (1997).

J.\,H.\,C.\,Whitehead:\,\,On certain sets of elements in a free group.
 Proc.\,\,London  Math.\,\,Soc.\,\,(2) \textbf{41}, 48--56  (1936-01).

 J.\,H.\,C.\,Whitehead:\,\,On equivalent sets of elements in a free group. Ann.\,\,of Math.\,\,(2)
\textbf{37}, 782--800  (1936-10).

Henry\,\,Wilton:\,\,Essential surfaces in graph pairs. J.\,\,Amer.\,\,Math.\,\,Soc.\,\,\textbf{31}, 893--919 (2018).

Ying-Qing\,\,Wu:\,\,Incompressible surfaces and Dehn surgery on $1$-bridge knots in handlebodies.
Math.\,\,Proc. Cambridge Philos.\,\,Soc.\,\,\textbf{120}, 687--696 (1996).

\leftskip 0pt \parindent 0pt

\end{document}